\documentclass[a4paper,12pt]{article}

\usepackage[utf8]{inputenc}
\usepackage{lmodern}
\usepackage[margin=2cm]{geometry}
\usepackage[english]{babel}
\usepackage{amsmath,amsthm,amssymb}
\usepackage{url}
\usepackage{hyperref}
\usepackage{graphicx}
\usepackage{caption}
\usepackage{subcaption}

\newtheorem{theorem}{Theorem}

\newtheorem{remark}[theorem]{Remark}

\newtheorem{lemma}[theorem]{Lemma}
\newtheorem{corollary}[theorem]{Corollary}

\def\cD{\mathcal{D}}

\def\cQ{\mathcal{Q}}

\def\Aut{\operatorname{Aut}}
\def\Vol{\operatorname{Vol}}
\def\cyl{\operatorname{cyl}}
\def\BigOh{\operatorname{O}}

\def\CC{\mathbb{C}}

\title{Asymptotics of lieanders with fixed composition sizes}
\author{Vincent Delecroix - (CNRS, LaBRI, Max-Planck Institut)}
\date{December 2018}

\begin{document}
\maketitle

\abstract{Lieanders are special cases of meanders and first appeared
in connection with Lie algebras. Using the results
from~\cite{DelecroixGoujardZografZorich} we prove a polynomial asymptotics
as $n \to \infty$ for the sequences $\left(L^{k^+,k^-}_n\right)_n$ counting
lieanders with $n$ arches and compositions of sizes $k^+$ and $k^-$.}

\section{Lieanders}
Meanders are combinatorial configurations of pairs of curves on the sphere.
According to~\cite{Lando:Zvonkin} the notion of meander was suggested by
V.~I.~Arnold~\cite{Arnold} though meanders were studied
already by H.~Poincar\'e~\cite{Poincare}. They are an object of particular
interest in statistical physics and the main conjecture concerning
their asymptotics is widely open~\cite{FGG1,FGG2,FGG3,FGG4}.
We are interested in a special kind of meanders that appeared in the context
of Lie algebras (see Degarchev and Kirillov~\cite{DegarchevKirillov} based
on some unpublished work of A.~G.~Elashvili and M.~Jibladze~\cite{ElashviliJibladze})
that we call lieanders\footnote{The name "lieander" comes
from~\cite{ElashviliJibladze}. Though what we call lieander in this note
would be a lieander of index 1 in their terminology. A restricted class
of lieanders are also called "birainbow meander" in~\cite{KarnauhovaLiebscher}.}.

A \emph{composition} of a non-negative integer $n$ is a sequence
of positive integers with sum $n$. The length of the composition
will be denoted $k$. The number of compositions of $n$ of length $k$
is $\binom{n-1}{k-1}$.

\begin{figure}[!ht]
\begin{center}\includegraphics{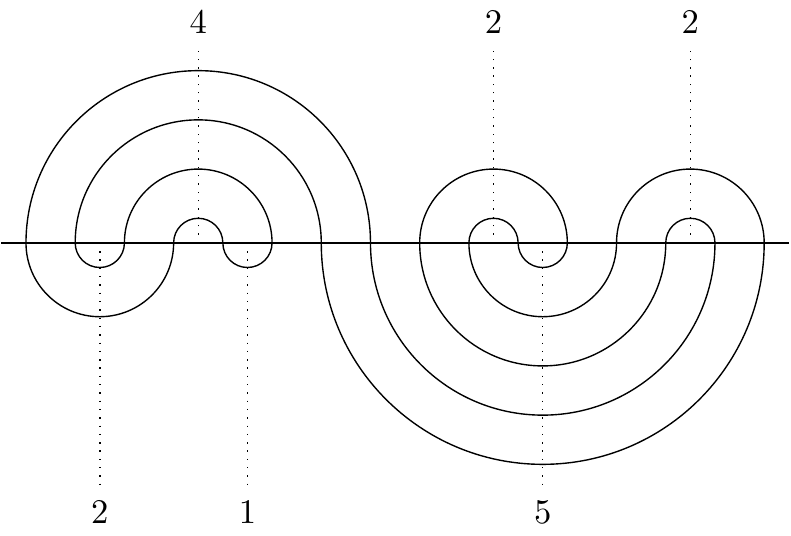}\end{center}
\caption{The $(3,3)$-lieander defined by the pair of compositions $((4,2,2),(2,1,5))$.}
\label{fig:lieander}
\end{figure}

Given two compositions $(c^+, c^-)$ of the same integer $n$ we
consider an upward sequence of nested arches given by the
composition $c^+$ and a downward sequence of nested arches
given by the composition $c^-$. Gluing together the pair of
configurations of arches we obtain a multicurve crossing the
horizontal line $2n$ times. See Figure~\ref{fig:lieander}.

\begin{figure}[!ht]
\begin{center}\includegraphics{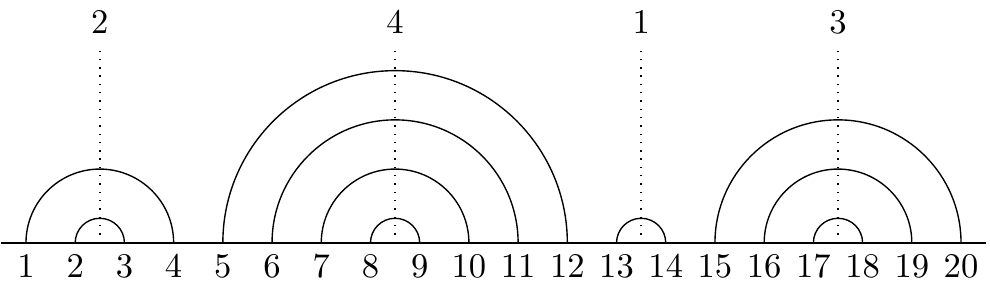}\end{center}
\caption{The arch configuration associated to the composition $(2,4,1,3)$.}
\label{fig:involution}
\end{figure}

More formally, to a composition $c = (c_1, c_2, \ldots, c_k)$ of $n$ we associate an
involution without fixed point $\sigma_c$ on $\{1,2,\ldots,2n\}$. We first split
the segment $\{1,2,\ldots,2n\}$ into pieces of size $2 \cdot c_1$, $2 \cdot c_2$, \ldots,
$2 \cdot c_k$ as follows
\begin{align*}
I_1 &:= \{1, 2, \ldots, 2 \cdot c_1 \}, \\
I_2 &:= \{2 \cdot c_1 + 1, \ldots, 2 \cdot (c_1 + c_2) \}, \\
\ldots \\
I_k &:= \{2 \cdot (c_1 + \ldots c_{k-1}) + 1, \ldots, 2 \cdot (c_1 + \ldots + c_k)\}.
\end{align*}
Then $\sigma_c$ is defined as the involution without fixed point reversing the
order on each $I_j$. To construct the $n$ arches from $\sigma_c$, simply connect $i$ to
$\sigma_c(i)$. See Figure~\ref{fig:involution}.

The connected component of the multicurve associated to a pair
$(c^+, c^-)$ are the orbits of $\langle \sigma_{c^+}, \sigma_{c^-} \rangle$.
When the multicurve associated to $(c^+, c^-)$ is connected, or equivalently when
the composition $\sigma_{c^+} \sigma_{c^-}$ is a product of two $n$-cycles,
we say that $(c^+, c^-)$ is a \emph{lieander}. A $(k^+, k^-)$-lieander is a
lieander $(c_+, c_-)$
made of compositions of lengths respectively $k^+$ and $k^-$.
We denote by $L_n^{k^+,k^-}$ the number of $(k^+,k^-)$-lieanders with sum $n$.
The main object of this note is to prove the following result.
\begin{theorem}
\label{thm:lieander:asymptotics}
Let $k^+$ and $k^-$ be two positive integers not both equal to $1$.
Then the number of $(k^+,k^-)$-lieanders with sum not greater than $x$ satisfies
the following asymptotics as $x$ tends to $\infty$
\begin{equation} \label{eq:asymptotics}
\sum_{n \leq x} L_n^{k^+,k^-} \sim
\delta_1 (k^+ + k^- - 4) \cdot \frac{x^{k^++k^--1}}{(k^++k^--1)\,(k^+-1)!\,(k^--1)!},
\end{equation}
where $\delta_1(k)$ is the real number
\[
\delta_1(k) = (k+4)\ (k+3) \cdot
\frac{(k+1)!!}{k!!} \cdot
\frac{1}{\pi^{k+2}} \cdot
\begin{cases}
\frac{1}{\pi}  & \text{if $k$ is odd}\\
\frac{1}{2}  &\text{if $k$ is even}
\end{cases}
\]
where $k!!$ is the double factorial
\[
(2n)!! = (2n)(2n-2)\cdots(2)
\qquad \text{and} \qquad
(2n+1)!! = (2n+1)(2n-1)\cdots(1).
\]
\end{theorem}
The number of pairs of compositions of size $(k^+, k^-)$ and sum not greater than
$x$ and of sizes $(k^+, k^-)$ has asymptotics the second term of the right-hand side
of~\eqref{eq:asymptotics}. Namely
\[
\frac{x^{k^++k^--1}}{(k^++k^--1)\,(k^+-1)!\,(k^--1)!}.
\]
In other words, $\delta_1(k^+ + k^- - 4)$ is the asymptotic density of the
number of lieanders among the pairs of compositions of size $(k^+, k^-)$.
The proof is an application of the work of the author with E. Goujard, P. Zograf and
A. Zorich~\cite{DelecroixGoujardZografZorich} and the number $\delta_1(k)$
is related to one-cylinder square-tiled surfaces in the minimal
strata of quadratic differential on the sphere $\cQ(k, -1^{k+4})$.

The constant $\delta_1(k)$ is easily evaluated
\[
\delta_1(-1) = \delta_1(0) = \frac{6}{\pi^2} \simeq 0.6079,\quad
\delta_1(1) = \frac{40}{\pi^4} \simeq 0.4106,\quad
\delta_1(2) = \frac{42}{2 \pi^4} \simeq 0.2310
\]
and applying Stirling formula, it is easy to prove that
\[
\frac{(k+1)!!}{k!!}
\cdot \begin{cases}
\frac{1}{\pi}  & \text{if $k$ is odd}\\
\frac{1}{2}  &\text{if $k$ is even}
\end{cases}
\sim
\sqrt{\frac{k}{2 \pi}}.
\]
Hence, as $k \to \infty$, we have
\[
\delta_1(k) \sim \frac{1}{\sqrt{2 \pi^5}} \cdot \left( \frac{1}{\pi}\right)^{k} \cdot k^{5/2} 
\]

\begin{remark}
For $(k,1)$-lieanders M. Duflo and R. Yu~\cite{DufloYu}
proved the existence of an exact polynomial behavior in a
different regime. More precisely, for each positive integer $n$, the sequence
$\left(L^{k,1}_{n+2k}\right)_k$ coincide with a polynomial for $k$ large enough.
\end{remark}

\noindent \textbf{Acknowledgment:}
The author wishes to thank A.~G.~Elashvili, E.~Goujard, P.~Hubert,
M.~Lapointe and A.~Zorich for fruitful discussions and remarks.

The numbers and graphics in Sections~\ref{sec:closed:lieandric:numbers}
and~\ref{sec:experimental}
were computed with the SageMath software~\cite{Sage} and the surface dynamics
package~\cite{SurfaceDynamics} of the author.

\section{Formulas and finer asymptotics for $k^+ + k^- \leq 4$}
In this section we describe precisely the pairs of compositions which give rise to
lieanders for the values $(k^+, k^-)$ with $k^+ + k^- \leq 4$. The condition
involves $\gcd$ which is natural as they correspond, up to a canonical
double cover, to discrete rotations. The corresponding counting formula involves
the Euler totient function $\phi$. Because in these situations the sequence
$(L_n^{k^+,k^-})_n$ is expicit, we obtain error terms in the asymptotics
\eqref{eq:asymptotics} of Theorem~\ref{thm:lieander:asymptotics}.

\begin{lemma}
\label{lem:description:rotation}
The following table provides necessary and sufficient conditions for
a $(k^+, k^-)$-pair of compositions of length $n$ to be a lieander
\begin{center}\begin{tabular}{lll}
$(k^+,k^-)$ & $(c^+, c^-)$ & condition \\ \hline
$(2,1)$ &
  $((a,b),(c))$ &
  $\gcd(a,n) = 1$ \\
$(3,1)$ &
 $((a,b,c),(d))$ &
 $\gcd(a+b,b+c) = 1$ \\
$(2,2)$ &
 $((a,b),(c,d))$ &
 $\gcd(a+d,n) = 1$.
\end{tabular}\end{center}
\end{lemma}
Such simple description of lieanders involving $\gcd$ do not exist for
other $k^+, k^-$, see~\cite{KarnauhovaLiebscher}.

\begin{proof}
These formulas already appeared in different places such as~\cite{KarnauhovaLiebscher}.
or~\cite{PakRedlich}. We only sketch a short unfied proof via canonical double covering.

As we have done in the introduction, the pair of compositions give rise to two
involutions without fixed points $\sigma^+$, $\sigma^-$ on $\{1, 2, \ldots, 2n\}$
that exchange even and odd numbers. The orbits of $\langle \sigma^+, \sigma^- \rangle$
are in bijections with the connected components of the multicurve built from the pair of
arches. These components are also in bijection with the cycles in the cyle decomposition
of the product $\sigma^+ \sigma^-$ restricted to even numbers.

For $(2,1)$-, $(3,1)$- and $(2,2)$- lieanders the product $\sigma^+ \sigma^-$ restricted
to even numbers correspond to an interval exchange transformations of rotation type
(i.e. that are either a rotation or an induced of rotations on 2 or 3 subintervals). The
lengths of the subintervals of the interval exchange are
easily expressed in terms of the parts of the compositions. Now, an interval exchange
of rotation type is made of a single orbit if and only if the underlying rotation 
is primitive, i.e. corresponds to the addition $x\mapsto x+m$ mod $\ell$ with $m$
relatively prime to $\ell$. From this remark it is easy to deduce the three formulas.
\end{proof}

\begin{corollary}
We have
\[
L_n^{2,1} = \phi(n)
\qquad
L_n^{3,1} = \sum_{k=1}^{n-1} \phi(k) - \phi(n),
\qquad \text{and} \qquad
L_n^{2,2} = (n-2) \phi(n)
\]
where $\phi(n)$ is the Euler totient function counting the number of
positive integers smaller than $n$ and relatively prime to $n$. 
\end{corollary}

Since the asymptotics of the Euler $\phi$ function is well known
(see e.g.~\cite[Section 18.5]{HardyWright}), we deduce the following corollary
on asymptotics.
\begin{corollary}
\label{cor:asymptotics:rotation}
As $x \to \infty$ we have the following asymptotics
\begin{align*}
\sum_{n \leq x} L_n^{2,1} &= \frac{6}{\pi^2} \cdot \frac{x^2}{2} + \BigOh(x\,\log(x)), \\
\sum_{n \leq x} L_n^{3,1} &= \frac{6}{\pi^2} \cdot \frac{x^3}{6} + \BigOh(x^2\,\log(x)), \\
\sum_{n \leq x} L_n^{2,2} &= \frac{6}{\pi^2} \cdot \frac{x^3}{3} + \BigOh(x^2\,\log(x)).
\end{align*}
\end{corollary}
The main terms of the asymptotics match our Theorem~\ref{thm:lieander:asymptotics} since
$\delta_1(-1) = \delta_1(0) = \frac{6}{\pi^2}$.

\section{Proof of the asymptotic in $n$ (Theorem~\ref{thm:lieander:asymptotics})}

We follow the definitions and notations from~\cite{DelecroixGoujardZografZorich}.

Let us fix the sizes $(k^+, k^-)$ for our pairs of compositions.
Contrarily to the count of meanders with fixed number of minimal
arches as in~\cite{DelecroixGoujardZografZorich}, the number of lieanders of
fixed sizes is not easily expressed in terms of square-tiled surfaces.
A pair of composition of sizes $(k^+, k^-)$ corresponds to a linear
involution with generalized permutation
\begin{equation}
\label{eq:generalized:perm}
\left(
\begin{array}{ll}
A_1\ A_1\ A_2\ A_2\ \ldots\ A_{k^+}\ A_{k^+} \\
B_1\ B_1\ B_2\ B_2\ \ldots\ B_{k^-}\ B_{k^-}
\end{array}
\right)
\end{equation}
and integral lengths given by the compositions. The stratum of
quadratic differentials of the generalized permutation~\eqref{eq:generalized:perm} is the
minimal stratum $\cQ(k, -1^{k+4})$ with $k=k^+ + k^- - 4$. The main result
of~\cite{DelecroixGoujardZografZorich} states that the asymptotic density
that such integral linear involution is made of a single orbit exists
and only depends on the stratum $\cQ(k, -1^{k+4})$. Moreover, the value of
this asymptotic density is $\delta_1(k)$ given by the formula
\begin{equation} \label{eq:delta1}
\delta_1(k) = \frac{\cyl_1(\cQ(k,-1^{k+4}))}{\Vol_1 (\cQ_1(k, -1^{k+4}))}
\end{equation}
where $\cyl_1(\cQ(k,-1^{k+4}))$ is the volume contribution of square-tiled surfaces
in the stratum $\cQ(k,-1^{k+4})$ with one horizontal cylinder of height $1/2$
and $\Vol_1 (\cQ_1(k,-1^{k+4}))$ is the Masur--Veech volume of the unit
hyperboloid $\cQ_1(k,-1^{k+4})$ in the stratum (or
equivalently the volume contribution of all square tiled surfaces).

From~\cite{AthreyaEskinZorich}, the volume of the minimal strata on the sphere is given by
\begin{equation} \label{eq:vol}
\Vol \cQ(k, -1^{k+4}) = 2 \pi^{k+2} \frac{k!!}{(k+1)!!}
\left\{
\begin{array}{ll}
\pi & \text{if $k$ is odd} \\
2   & \text{otherwise.}
\end{array} \right.
\end{equation}
Which is the denominator in~\eqref{eq:delta1}.

To compute the numerator in~\eqref{eq:delta1}, we consider surfaces in
$\cQ(k, -1^{k+4})$ which are made of a single horizontal
cylinder. This stratum is particularly simple as there is a single one
cylinder separatrix diagram that is obtained by putting two poles on one side,
and the zeros and the other poles on the other
side. The volume contribution can be computed in at least two ways. We can
first use the formula for a given separatrix diagram $\cD$ from~\cite{DelecroixGoujardZografZorich-Yoccoz}
\[
\cyl_1(\cD) =
\frac{2^{s+2}}{|\Aut(\cD)|} \cdot \frac{(t+u-2)!}{(t-1)!(u-1)!} \cdot
\frac{\mu_{-1}! \mu_1! \mu_2! \cdots}{(d-2)!}
\]
where $s$, $t$ and $u$ are respectively the number of saddle connections
on the bottom and top of $\cD$, twice on top and twice on bottom, $\cQ(-1^{\mu_{-1}},
0^{\mu_0}, 1^{\mu_1}, \ldots)$ is the stratum of $\cD$ and $d = 2g + \# \textrm zeros - 2$
is the dimension of the stratum.
Replacing $s = 0$, $t=m+2$, $u = 1$, $d = m+3$ and $\mu = (-1^{k+4}, k)$ in the formula
we obtain
\begin{equation} \label{eq:cyl1}
\cyl_1(\cD) = 2 (k+4) (k+3).
\end{equation}

An alternative computation is available in the very special case of a
stratum of genus 0 thanks to~\cite{DelecroixGoujardZografZorich}.
We have
\[
\cyl_1( \cQ(\mu) ) = 2 \sum_{\mu \subset \nu}
\binom{|\nu|+4}{|\mu|+2} \binom{\nu_0}{\mu_0} \binom{\nu_1}{\mu_1} \cdots
\]
In the special case of interest to us we obtain
\begin{equation} \label{eq:cyl1:alt}
\cyl_1(\cQ(k, -1^{k+4})) =
2 \left( \binom{k+4}{k+2} + \binom{k+4}{2} \right) = 2 (k+4) (k+3)
\end{equation}
which coincides with our previous formula.

Gathering~\eqref{eq:vol} and~\eqref{eq:cyl1} in Formula~\eqref{eq:delta1} we obtain
\[
\delta_1(k) = \frac{2 (k+4) (k+3)}%
{2 \pi^{k+2} \frac{k!!}{k+1!!} \left\{
\begin{array}{ll}
\pi & \text{if $k$ is odd} \\
2 & \text{otherwise}
\end{array} \right.
}.
\]

\newpage
\section{Some numbers for $k^+ + k^- \leq 7$}
\label{sec:closed:lieandric:numbers}
In the array below we gathered numbers $L_n^{k^+,k^-}$ for
$k^+ + k^- \leq 7$ and $n \leq 50$. \\ \bigskip
\begin{center}%
\resizebox{0.8\textwidth}{!}{%
\begin{tabular}{r|r|rr|rr|rrr|rrr|}
n  & $L_n^{2,1}$
   & $L_n^{3,1}$ & $L_n^{2,2}$
   & $L_n^{4,1}$ & $L_n^{3,2}$
   & $L_n^{5,1}$ & $L_n^{4,2}$ & $L_n^{3,3}$
   & $L_n^{6,1}$ & $L_n^{5,2}$ & $L_n^{4,3}$ \\ \hline
2 & 1 & 0 & 0 & 0 & 0 & 0 & 0 & 0 & 0 & 0 & 0\\
3 & 2 & 0 & 2 & 0 & 0 & 0 & 0 & 0 & 0 & 0 & 0\\
4 & 2 & 2 & 4 & 0 & 1 & 0 & 0 & 0 & 0 & 0 & 0\\
5 & 4 & 2 & 12 & 0 & 4 & 0 & 0 & 2 & 0 & 0 & 0\\
6 & 2 & 8 & 8 & 2 & 13 & 0 & 2 & 4 & 0 & 0 & 1\\
7 & 6 & 6 & 30 & 2 & 24 & 0 & 6 & 22 & 0 & 0 & 4\\
8 & 4 & 14 & 24 & 12 & 37 & 2 & 24 & 40 & 0 & 3 & 15\\
9 & 6 & 16 & 42 & 6 & 80 & 2 & 34 & 96 & 0 & 8 & 48\\
10 & 4 & 24 & 32 & 34 & 93 & 12 & 94 & 140 & 2 & 33 & 101\\
11 & 10 & 22 & 90 & 22 & 158 & 6 & 130 & 318 & 2 & 44 & 226\\
12 & 4 & 38 & 40 & 64 & 181 & 50 & 256 & 368 & 12 & 155 & 395\\
13 & 12 & 34 & 132 & 44 & 302 & 24 & 332 & 738 & 6 & 176 & 774\\
14 & 6 & 52 & 72 & 126 & 299 & 118 & 574 & 812 & 54 & 507 & 1125\\
15 & 8 & 56 & 104 & 72 & 544 & 62 & 704 & 1496 & 26 & 518 & 2032\\
16 & 8 & 64 & 112 & 200 & 473 & 250 & 1104 & 1552 & 170 & 1255 & 2715\\
17 & 16 & 64 & 240 & 132 & 782 & 120 & 1372 & 2832 & 72 & 1270 & 4542\\
18 & 6 & 90 & 96 & 302 & 687 & 488 & 1970 & 2652 & 396 & 2777 & 5673\\
19 & 18 & 84 & 306 & 202 & 1152 & 224 & 2366 & 4772 & 170 & 2646 & 9228\\
20 & 8 & 112 & 144 & 428 & 1037 & 820 & 3164 & 4372 & 860 & 5435 & 10797\\
21 & 12 & 116 & 228 & 268 & 1672 & 412 & 3860 & 7452 & 340 & 5128 & 17016\\
22 & 10 & 130 & 200 & 626 & 1361 & 1336 & 4986 & 6724 & 1648 & 9879 & 19029\\
23 & 22 & 128 & 462 & 394 & 2200 & 642 & 5950 & 11604 & 660 & 9218 & 29324\\
24 & 8 & 164 & 176 & 788 & 1783 & 2002 & 7372 & 9912 & 2978 & 16533 & 31759\\
25 & 20 & 160 & 460 & 536 & 2984 & 970 & 8760 & 16808 & 1174 & 15500 & 48088\\
26 & 12 & 188 & 288 & 1098 & 2341 & 3006 & 10582 & 14220 & 4948 & 27049 & 50341\\
27 & 18 & 194 & 450 & 678 & 3724 & 1454 & 12562 & 23462 & 1946 & 24646 & 75364\\
28 & 12 & 218 & 312 & 1352 & 3049 & 4168 & 14712 & 19776 & 7882 & 41217 & 77141\\
29 & 28 & 214 & 756 & 904 & 4742 & 2026 & 17448 & 32794 & 3112 & 38054 & 113610\\
30 & 8 & 262 & 224 & 1678 & 3593 & 5938 & 19938 & 26296 & 12178 & 62107 & 113371\\
31 & 30 & 248 & 870 & 1126 & 5926 & 2822 & 23546 & 43884 & 4756 & 56618 & 166526\\
32 & 16 & 292 & 480 & 2100 & 4619 & 7772 & 26524 & 35476 & 17868 & 89423 & 163477\\
33 & 20 & 304 & 620 & 1332 & 7136 & 3862 & 31092 & 56844 & 7006 & 81778 & 236796\\
34 & 16 & 328 & 512 & 2626 & 5537 & 10488 & 34534 & 45820 & 26042 & 126957 & 228599\\
35 & 24 & 336 & 792 & 1684 & 9006 & 5016 & 40388 & 74696 & 10140 & 114974 & 329734\\
36 & 12 & 372 & 408 & 3012 & 6465 & 13186 & 44084 & 58172 & 36318 & 173649 & 313349\\
37 & 36 & 360 & 1260 & 2056 & 10458 & 6508 & 52000 & 95152 & 14220 & 158612 & 449562\\
38 & 18 & 414 & 648 & 3758 & 7857 & 17330 & 56038 & 73972 & 50342 & 237595 & 421683\\
39 & 24 & 426 & 888 & 2356 & 12060 & 8350 & 65148 & 117714 & 19580 & 214200 & 600996\\
40 & 16 & 458 & 608 & 4320 & 9391 & 21258 & 69816 & 92472 & 67028 & 312597 & 559705\\
41 & 40 & 450 & 1560 & 2888 & 14508 & 10462 & 81592 & 148202 & 26348 & 284788 & 793208\\
42 & 12 & 518 & 480 & 4998 & 10375 & 27032 & 85074 & 111708 & 90216 & 412583 & 726155\\
43 & 42 & 500 & 1722 & 3374 & 16924 & 13028 & 100378 & 181536 & 34930 & 372092 & 1030880\\
44 & 20 & 564 & 840 & 5888 & 12675 & 32454 & 105288 & 138840 & 115926 & 530871 & 938347\\
45 & 24 & 580 & 1032 & 3708 & 19016 & 16122 & 121308 & 216664 & 45706 & 480528 & 1317228\\
46 & 22 & 606 & 968 & 6930 & 14337 & 40360 & 126866 & 166660 & 151900 & 680025 & 1190891\\
47 & 46 & 604 & 2070 & 4510 & 22382 & 19450 & 147706 & 265940 & 58876 & 612920 & 1673598\\
48 & 16 & 680 & 736 & 7512 & 15829 & 47340 & 150328 & 197636 & 190596 & 851655 & 1495615\\
49 & 42 & 670 & 1974 & 5178 & 25826 & 23480 & 176774 & 316646 & 74988 & 773180 & 2097602\\
50 & 20 & 734 & 960 & 8926 & 18661 & 58238 & 180594 & 236660 & 243210 & 1070837 & 1862765\\

\end{tabular}%
}%
\end{center}

\newpage
\section{Empirical asymptotics for $k^+ + k^- = 5$}
\label{sec:experimental}
In this section we report on experiments regarding the squences
$(L_n^{k^+,k^-})_n$ with $k^+ + k^- = 5$, that is $(k^+, k^-) \in \{(4,1), (3,2)\}$.
These sizes correspond to the stratum $\cQ(1,-1^5)$ and the two experiments
match our asymptotics from Theorem~\ref{thm:lieander:asymptotics}.

\begin{figure}[!ht]
\begin{subfigure}[c]{0.5\textwidth}
\includegraphics{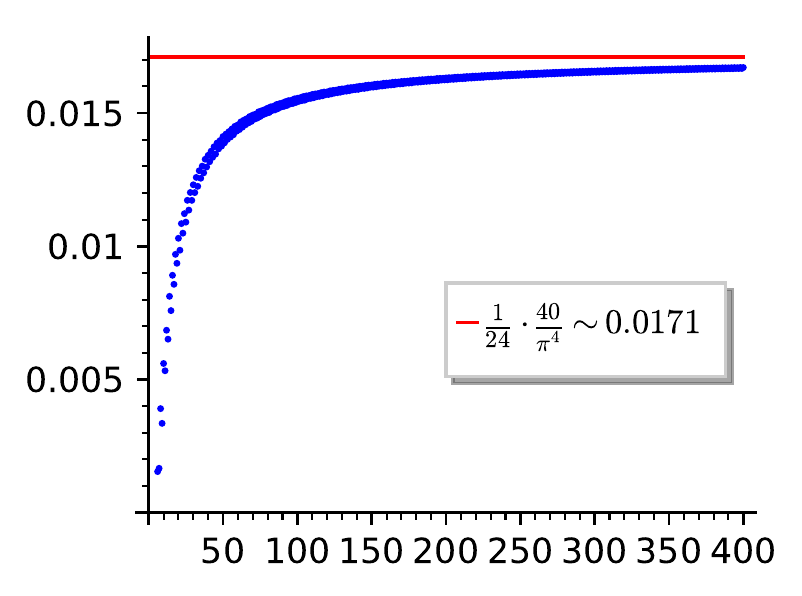}
\subcaption{$(k^+, k^-) = (4,1)$}
\end{subfigure}%
\begin{subfigure}[c]{0.5\textwidth}
\includegraphics{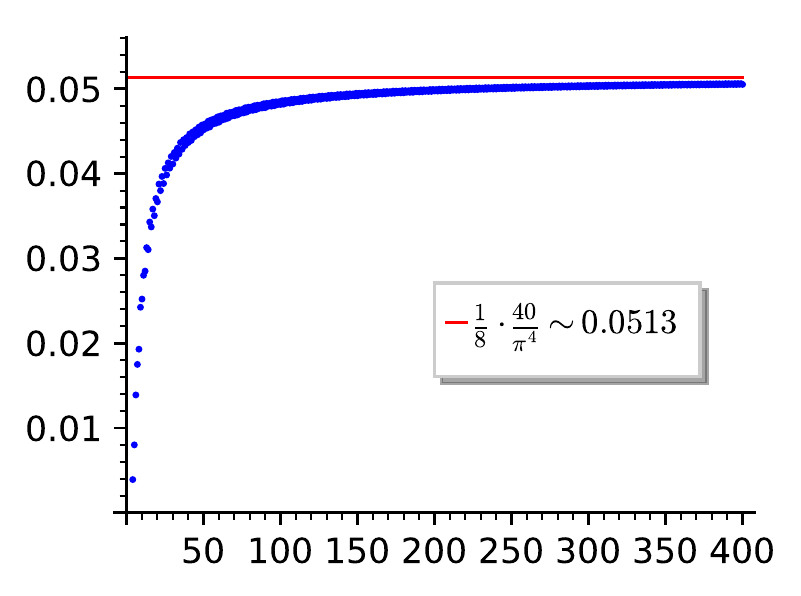}
\subcaption{$(k^+, k^-) = (3,2)$}
\end{subfigure}%
\caption{Graphics of $x \mapsto \displaystyle \frac{1}{x^4} \left(\sum_{n \leq x} L_n^{(k^+,k^-)}\right)$ for $k^+ + k^- = 5$ (blue points). We also indicate the two limit values that are rational multiples of $\delta_1(1) = \frac{40}{\pi^4}$ (red line).}
\label{fig:graphics:cum}
\end{figure}

\begin{figure}[!ht]
\begin{subfigure}[c]{0.5\textwidth}
\includegraphics{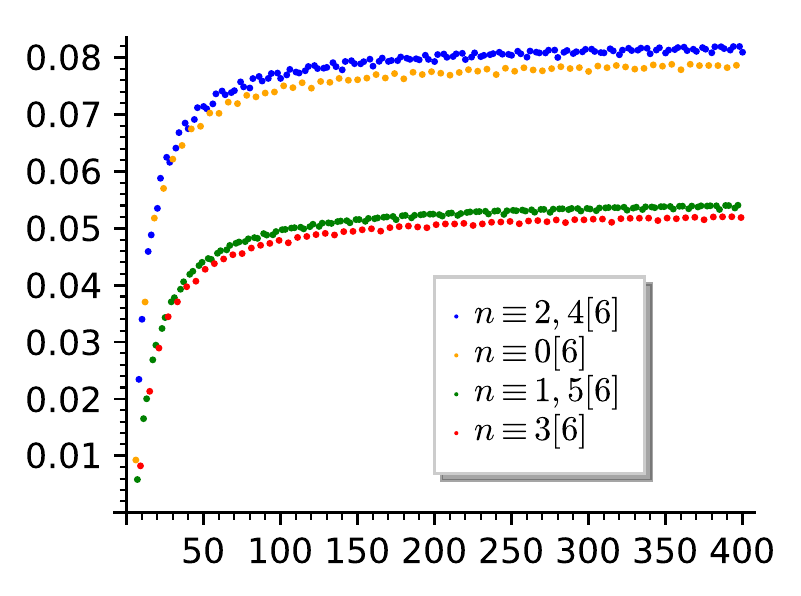}
\subcaption{$(k^+, k^-) = (4,1)$}
\end{subfigure}%
\begin{subfigure}[c]{0.5\textwidth}
\includegraphics{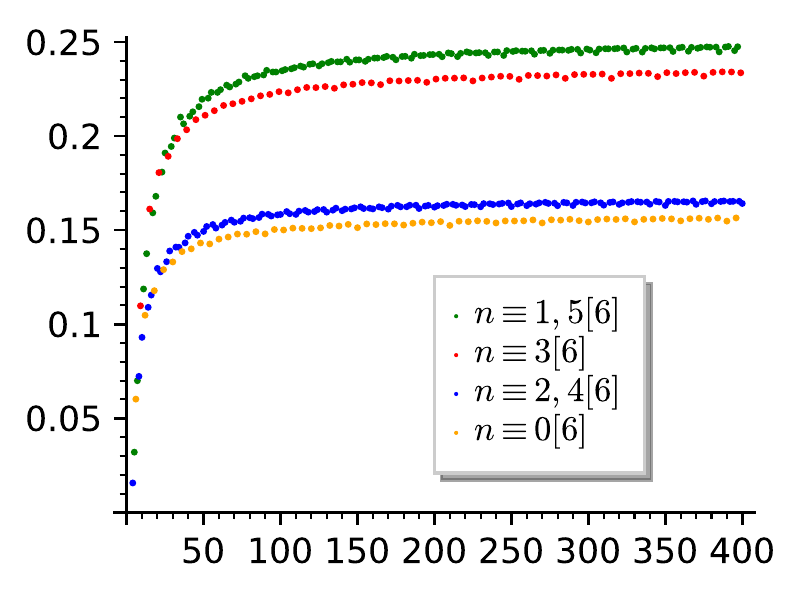}
\subcaption{$(k^+, k^-) = (3,2)$}
\end{subfigure}%
\caption{Graphics of $L_n^{(k^+,k^-)}$ with 4 colours depending on the residue modulo
6 of $n$: orange for $0$, green for $1$ and $6$, blue for $2$ and $4$, red for $3$.}
\label{fig:graphics:mod}
\end{figure}

Let us make three remarks about Figure~\ref{fig:graphics:cum} that are not explained
by our proof of Theorem~\ref{thm:lieander:asymptotics}. First of all it seems
that the convergence to the limit is from below. Secondly, one
can notice some oscillations which indicates some sensitivity
of $L_n^{(k^+, k^-)}$ depending on the factorization of $n$. Figure~\ref{fig:graphics:mod}
is a graphic of the sequence $L_n^{(k^+, k^-)} / n^3$ (no summation) that tend to confirm
this phenomenon. Though, the relative positions of the congruence classes are not the
same for $(4,1)$-lieanders and $(3,2)$-lieanders. Finally,
we have no clue about the speed of convergence of
$\frac{1}{x^4} \sum_{n \leq x} L_n^{(k^+, k^-)}$. Figure~\ref{fig:graphics:mod}
indicates that $L_n^{(k^+,k^-)} / n^3$ is rather mild. One can hope
as in Corollary~\ref{cor:asymptotics:rotation} that the convergence speed
is $\log(x)/x$.

\end{document}